\documentclass[fleqn]{mat01}
\usepackage{times,mathtimy,amssymb,latexsym,amscd}
\begin{document}

\setcounter{page}{499} \firstpage{499}

%\input{tcilatex}
%\theoremstyle{definition}
%\theoremstyle{remark}
%\numberwithin{equation}{section}
%\newcommand{\thmref}[1]{Theorem~\ref{#1}}
%\newcommand{\secref}[1]{\S\ref{#1}}
%\newcommand{\lemref}[1]{Lemma~\ref{#1}}
%\input tcilatex

\newtheorem{theoree}{\bf Theorem}
\renewcommand\thetheoree{\Alph{theoree}}

\newtheorem{theore}{Theorem}
\renewcommand\thetheore{\arabic{section}.\arabic{theore}}
\newtheorem{theor}[theore]{\bf Theorem}
\newtheorem{lem}[theore]{\it Lemma}
\newtheorem{propo}[theore]{\rm PROPOSITION}
\newtheorem{coro}[theore]{\rm COROLLARY}
\newtheorem{definit}[theore]{\rm DEFINITION}
\newtheorem{probl}[theore]{\it Problem}
\newtheorem{exampl}[theore]{\it Example}

\def\remar{\trivlist \item[\hskip \labelsep{\it Remark.}]}

\renewcommand\theequation{\arabic{section}.\arabic{equation}}

\title{$\pmb{A}$-Statistical extension of the Korovkin type approximation theorem}

\markboth{Esra Erku\c{s} and Oktay Duman}{$A$-Statistical
extension of the Korovkin theorem}

\author{ESRA ERKU\c{S}$^{1}$ and OKTAY DUMAN$^{2}$}

\address{$^{1}$Department of Mathematics, Faculty of Sciences and
Arts, \c{C}anakkale Onsekiz Mart University, Terzio\u{g}lu
Kamp\"{u}s\"{u}~17020, \c{C}anakkale, Turkey\\
\noindent $^{2}$Faculty of Arts and Sciences, Department of
Mathematics, TOBB University of Economy and Technology,
S\"{o}\u{g}\"{u}t\"{o}z\"{u} 06530, Ankara, Turkey\\
\noindent E-mail: erkus@comu.edu.tr; oduman@etu.edu.tr}

\volume{115}

\mon{November}

\parts{4}

\pubyear{2005}

\Date{MS received 3 October 2004}

\begin{abstract}
In this paper, using the concept of $A$-statistical convergence
which is a regular (non-matrix) summability method, we obtain a
general Korovkin type approximation theorem which concerns the
problem of approximating a function $f$ by means of a sequence
$\{L_{n}f\}$ of positive linear operators.
\end{abstract}


\maketitle

\section{Introduction}

Let $A=(a_{jn})$ be an infinite summability matrix. For a given
sequence $x:=(x_{n}),$ the $A$-{\it transform of} $x$, denoted by
$Ax:=((Ax)_{j})$, is given by $(Ax)_{j}=\sum_{n=1}^{\infty}
a_{jn}x_{n}$ provided the series converges for each $j\in
\Bbb{N}$, the set of all natural numbers. We say that $A$ is {\it
regular} if $\lim Ax=L$ whenever $\lim x=L$ \cite{Boos}. Assume
that $A$ is a non-negative regular summability matrix. Then
$x=(x_{n})$ is said to be $A$-{\it statistically convergent to}
$L$ if, for every $\varepsilon >0$, $\lim_{j}\sum_{n:\left|
x_{n}-L\right| \geq \varepsilon}a_{jn}=0,$ which is denoted by
$st_{A}-\lim x=L$ \cite{Fre} (see also \cite{Kol,Mil}). We note
that by taking $ A=C_{1}$, the Ces\`{a}ro matrix, $A$-statistical
convergence reduces to the concept of {\it statistical
convergence} (see \cite{Fas,Fri,Ste} for details). If $A$ is the
identity matrix, then $A$-statistical convergence coincides with
the ordinary convergence. It is not hard to see that every
convergent sequence is $A$-statistically convergent. However, Kolk
\cite{Kol} showed that $A$-statistical convergence is stronger
than convergence when $A=(a_{jn})$ is a regular summability matrix
such that $\lim_{j}\max_{n}| a_{jn}| =0.$ It should be noted that
$A$-statistical convergence may also be given in normed spaces
\cite{Kol-2}.

Approximation theory has important applications in the theory of
polynomial approximation, in various areas of functional analysis
\cite{Alt,Boj,Che,Khan,Kor}. The study of the Korovkin type
approximation theory is a well-established area of research, which
deals with the problem of approximating a function $f$ by means of
a sequence $\{ L_{n}f\} $ of positive linear operators.
Statistical convergence, which was introduced nearly fifty years
ago, has only recently become an area of active research.
Especially it has made an appearance in approximation theory (see,
for instance, \cite{Dum,Gad}). The aim of the present paper is to
investigate their use in approximation theory settings.

Throughout this paper $I:=[0,\infty)$. As usual, let
$C(I):=\{f\hbox{:}\ f$ is a real-valued continuous functions on $I
\},$ and $C_{B}(I):=\{f\in C(I)\hbox{:}\ f$ is bounded on $I \}.$
Consider the space $H_{w}$ of all real-valued functions $f$
defined on $I$ and satisfying
\begin{equation*}
\left| f(x)-f(y)\right| \leq w\left(f;\left|
\frac{x}{1+x}-\frac{y}{1+y} \right| \right),
\end{equation*}
where $w$ is the modulus of continuity given by, for any $\delta
>0,$
\begin{equation*}
w(f;\delta):=\sup_{\underset{\left| x-y\right| \leq \delta
}{x,y\in I} }\left| f(x)-f(y)\right|.
\end{equation*}
Assume that $L$ is a linear operator mapping $H_{w}$ into
$C_{B}(I)$. As usual, we say that $L$ is a {\it positive} linear
operator provided that $f\geq 0$ implies $Lf\geq 0, f\in H_{w}$.
Also, we denote the value of $ Lf $ at a point $x\in I$ by
$L(f(u);x)$ or simply $L(f;x)$.

We now recall that \c{C}akar and Gadjiev \cite{C-G} obtained the
following result.

\begin{theoree}[\!] Let $\{ L_{n}\}$ be a
sequence of positive linear operators from $H_{w}$ into
$C_{B}(I)$. Then{\rm ,} the uniform convergence of the sequence
$\{L_{n}f\}$ to $f$ on $I$ holds for any function $f\in H_{w}$ if
the sequences $\{ L_{n}\varphi_{i}\}$ converge uniformly to
$\varphi _{i},$ where $\varphi_{i}(u)=
\big(\frac{u}{1+u}\big)^{i},$ $i=0,1,2.$
\end{theoree}

\section{$\pmb{A}$-Statistical approximation}

In this section, replacing the limit operator by the
$A$-statistical limit operator and considering a sequence of
positive linear operators defined on the space of all real-valued
continuous and bounded functions on a subset $I^{m}$ of
$\Bbb{R}^{m}$, the real $m$-dimensional space, where
$I^{m}:=I\times I\times \cdots \times I,$ we give an extension of
Theorem~A.

To achieve this we first consider the case of $m=2$.

Let $K:=I^{2}=[0,\infty)\times [0,\infty).$ Then, the sup norm on
$ C_{B}(K) $ is given by
\begin{equation*}
\left\| f\right\| :=\sup\limits_{(x,y)\in K}\left| f(x,y)\right|,
\quad f\in C_{B}(K),
\end{equation*}
and also the value of $Lf$ at a point $(x,y)\in K$ is denoted by
$L(f(u,v);x,y)$ or simply $L(f;x,y)$.

We consider the modulus of continuity
$w_{2}(f;\delta_{1},\delta_{2})$ (for the functions of two
variables) given by, for any $\delta _{1},$ $\delta _{2}>0,$
\begin{align*}
w_{2}(f;\delta _{1},\delta _{2}) &:=\sup \{\left|
f(u,v)-f(x,y)\right|\hbox{:}\ (u,v), (x,y)\in K,\hbox{ and}\\[.2pc]
&\quad\, \left| u-x\right| \leq \delta _{1}, \left| v-x\right|
\leq \delta _{2}\}.
\end{align*}
It is clear that a necessary and sufficient condition for a
function $f\in C_{B}(K)$ is
\begin{equation*}
\lim_{\delta _{1}\rightarrow 0,\text{ }\delta _{2}\rightarrow
0}w_{2}(f;\delta _{1},\delta _{2})=0.
\end{equation*}
We now introduce the space $H_{w_{2}}$ of all real-valued
functions $f$ defined on $K$ and satisfying
\begin{equation}
\left| f(u,v)-f(x,y)\right| \leq w_{2}\left(f;\left|
\frac{u}{1+u}-\frac{x}{ 1+x}\right|,\left|
\frac{v}{1+v}-\frac{y}{1+y}\right| \right).
\end{equation}
Then observe that any function in $H_{w_{2}}$ is continuous and
bounded on $K.$ \pagebreak

With this terminology we have the following:

\begin{theor}[\!]
Let $A=(a_{jn})$ be a non-negative regular summability matrix{\rm
,} and let $\{ L_{n}\} $ be a sequence of positive linear
operators from $ H_{w_{2}}$ into $C_{B}(K)$. Then{\rm,} for any
$f\in H_{w_{2}},$
\begin{equation}
st_{A}-\lim\limits_{n}\left\| L_{n}f-f\right\| =0
\end{equation}
is satisfied if the following holds{\rm :}
\begin{equation}
st_{A}-\lim\limits_{n}\left\| L_{n}f_{i}-f_{i}\right\| =0,\quad
i=0,1,2,3,
\end{equation}
where
\begin{align*}
f_{0}(u,v) &= 1,\text{ }f_{1}(u,v)=\frac{u}{1+u},\text{
}f_{2}(u,v)=\frac{v}{1+v}, \\[.2pc]
f_{3}(u,v) &= \left(\frac{u}{1+u}\right)
^{2}+\left(\frac{v}{1+v}\right) ^{2}.
\end{align*}
\end{theor}

\begin{proof}
Assume that $(2.3)$ holds, and let $f\in H_{w_{2}}.$ By $(2.3),$
for every $ \varepsilon >0,$ there exist $\delta _{1},\delta
_{2}>0$ such that $| f(u,v)-f(x,y)| <\varepsilon $ holds for all
$(u,v)\in K$ satisfying $\big|\frac{u}{1+u}-\frac{x}{1+x}\big|
<\delta_{1}$ and $\big|\frac{v}{1+v} -\frac{y}{1+y}\big|
<\delta_{2}$. Let
\begin{align*}
K_{\delta _{1},\delta _{2}}:=\left\{ (u,v)\in K\hbox{:}\ \left|
\frac{u}{1+u}\!-\!\frac{x }{1+x}\right| <\delta _{1}\text{ and }
\left| \frac{v}{1+v}\!-\!\frac{y}{1+y} \right| <\delta
_{2}\right\}.
\end{align*}
Hence we may write
\begin{align}
\left| f(u,v)-f(x,y)\right| & =\left| f(u,v)-f(x,y)\right| \chi
_{K_{\delta_{1},\delta_{2} }}(u,v)\nonumber\\[.2pc]
&\quad\, +\left| f(u,v)-f(x,y)\right| \chi _{K\setminus
K_{\delta_{1},\delta_{2}
}}(u,v)\nonumber\\[.2pc]
&<\varepsilon +2M\chi _{K\setminus K_{\delta_{1},\delta_{2}
}}(u,v),
\end{align}
where $\chi _{R}$ denotes the characteristic function of the set
$R$ and $ M:=\| f\|$. We also get that
\begin{equation}
\chi _{K\setminus K_{\delta_{1},\delta_{2} }}(u,v)\leq
\frac{1}{\delta _{1}^{2}}\left(\frac{u}{1+u}-\frac{x}{1+x}\right)
^{2}+\frac{1}{\delta
_{2}^{2}}\left(\frac{v}{1+v}-\frac{y}{1+y}\right) ^{2}.
\end{equation}
Combining $(2.4)$ with $(2.5)$ we have
\begin{align}
\left| f(u,v)-f(x,y)\right| \leq \varepsilon \!+\!\frac{2M}{\delta
^{2}}\left\{ \left(\frac{u}{1+u}\!-\!\frac{x}{1+x}\right)
^{2}\!+\!\left(\frac{v}{1+v}\!-\!\frac{y}{
1+y}\right)^{2}\right\},
\end{align}
where $\delta :=\min \{\delta _{1},\delta _{2}\}.$ Using linearity
and positivity of the operators $L_{n}$ we get, for any $n\in
\Bbb{N}$,
\begin{align*}
\left| L_{n}(f;x,y)-f(x,y)\right| &\leq L_{n}\left(\left|
f(u,v)-f(x,y)\right| ;x,y\right)\\[.2pc]
&\quad\, +\left| f(x,y)\right| \left|
L_{n}(f_{0};x,y)-f_{0}(x,y)\right|.
\end{align*}
Then, by $(2.5),$ we obtain
\begin{align*}
\left| L_{n}(f;x,y)-f(x,y)\right| &\leq \varepsilon
L_{n}(f_{0};x,y)\\[.2pc]
&\quad\, +\frac{2M}{\delta ^{2}}\left\{ L_{n}\bigg(\left(\frac{u}{1+u}-\frac{x}{1+x}\right) ^{2};x,y\bigg) \right.\\[.2pc]
&\quad\, \left. +L_{n}\bigg(\left(\frac{v}{1+v}-
\frac{y}{1+y}\right)^{2};x,y\bigg) \right\} \\[.2pc]
&\quad\, +M\left| L_{n}(f_{0};x,y)-f_{0}(x,y)\right|.
\end{align*}
By some simple calculations we have
\begin{align*}
\left| L_{n}(f;x,y)-f(x,y)\right| &\leq \varepsilon +(\varepsilon
+M)\left| L_{n}(f_{0};x,y)-f_{0}(x,y)\right|\\[.2pc]
&\quad\, +\frac{2M}{\delta ^{2}}
\left\{L_{n}(f_{3};x,y)-\frac{2x}{1+x}
L_{n}(f_{1};x,y)\right. \\[.2pc]
&\quad\, -\frac{2y}{1+y}L_{n}(f_{2};x,y) \\[.2pc]
&\quad\, +\bigg(\left(\frac{x}{1+x}\right) ^{2}+\left(\frac{y}{1+y}\right) ^{2}\bigg) L_{n}(f_{0};x,y)\bigg\} \\[.2pc]
&= \varepsilon +(\varepsilon +M)\left| L_{n}(f_{0};x,y)-f_{0}(x,y)\right| \\[.2pc]
&\quad\, +\frac{2M}{\delta ^{2}} \left(L_{n}(f_{3};x,y)-f_{3}(x,y)\right) \\[.2pc]
&\quad\, -\frac{4M}{\delta ^{2}}\left(\frac{x}{1+x}\right) \left(L_{n}(f_{1};x,y)-f_{1}(x,y)\right) \\[.2pc]
&\quad\, -\frac{4M}{\delta ^{2}}\left(\frac{y}{1+y}\right) \left(L_{n}(f_{2};x,y)-f_{2}(x,y)\right) \\[.2pc]
&\quad\, +\frac{2M}{\delta ^{2}}\bigg(\left(\frac{x}{1+x}\right)
^{2}+\left(\frac{y}{1+y}\right) ^{2}\bigg) \\[.2pc]
&\quad\, \times \left(L_{n}(f_{0};x,y)-f_{0}(x,y)\right) \\[.2pc]
&\leq \varepsilon +\left(\varepsilon +M+\frac{4M}{\delta
^{2}}\right)\left|L_{n}(f_{0};x,y)- f_{0}(x,y)\right| \\[.2pc]
&\quad\, +\frac{4M}{\delta ^{2}}\left| L_{n}(f_{1};x,y)-f_{1}(x,y)\right| \\[.2pc]
&\quad\, +\frac{4M}{\delta ^{2}}\left| L_{n}(f_{2};x,y)-f_{2}(x,y)\right| \\[.2pc]
&\quad\, +\frac{2M}{\delta ^{2}}\left|
L_{n}(f_{3};x,y)-f_{3}(x,y)\right|.
\end{align*}
Then, taking supremum over $(x,y)\in K$ we have
\begin{align}
\left\| L_{n}f-f\right\| & \leq \varepsilon +B \{\left\|
L_{n}f_{0}-f_{0}\right\| +\left\| L_{n}f_{1}-f_{1}\right\|\nonumber\\[.2pc]
&\quad\, +\left\| L_{n}f_{2}-f_{2}\right\| +\left\|
L_{n}f_{3}-f_{3}\right\|\},
\end{align}
where $B:=\varepsilon +M+\frac{4M}{\delta ^{2}}.$ For a given
$r>0$, choose $\varepsilon >0$ such that $\varepsilon <r$. Define
the following sets:
\begin{align*}
D &:=\left\{ n\hbox{:}\ \left\| L_{n}f-f\right\| \geq r\right\}, \\[.2pc]
D_{1} &:=\left\{ n\hbox{:}\ \left\| L_{n}f_{0}-f_{0}\right\| \geq
\frac{r-\varepsilon }{4B}\right\}, \\[.2pc]
D_{2} &:=\left\{ n\hbox{:}\ \left\| L_{n}f_{1}-f_{1}\right\| \geq
\frac{r-\varepsilon }{4B}\right\}, \\[.2pc]
D_{3} &:=\left\{ n\hbox{:}\ \left\| L_{n}f_{2}-f_{2}\right\| \geq
\frac{r-\varepsilon }{4B}\right\}, \\[.2pc]
D_{4} &:=\left\{ n\hbox{:}\ \left\| L_{n}f_{3}-f_{3}\right\| \geq
\frac{ r-\varepsilon }{4B}\right\}.
\end{align*}
It follows from $(2.7)$ that $D\subseteq D_{1}\cup D_{2}\cup
D_{3}\cup D_{4}. $ Therefore, for each $j\in \Bbb{N}$, we may
write
\begin{equation}
\sum_{n\in D}a_{jn}\leq \sum_{n\in D_{1}}a_{jn}+\sum_{n\in
D_{2}}a_{jn}+\sum_{n\in D_{3}}a_{jn}+\sum_{n\in D_{4}}a_{jn}.
\end{equation}
Letting $j\rightarrow \infty $ in (2.8) and using (2.3) we
conclude that
\begin{equation*}
\lim_{j}\sum_{n\in D}a_{jn}=0,
\end{equation*}
whence gives (2.2). So the proof is completed.\hfill $\Box$
\end{proof}

Now replace $I^{2}$ by $I^{m}:=[0,\infty)\times \cdots \times
[0,\infty)$ and consider the modulus of continuity $w_{m}(f;\delta
_{1},\dots,\delta _{m})$ (for the functions $f$ of $m$-variables)
given by, for any $\delta _{1}, \dots,\delta _{m}>0,$
\begin{align*}
w_{m}(f;\delta _{1},\dots,\delta _{m}) &:=\sup \{|
f(u_{1},\dots,u_{m})-f(x_{1},\dots,x_{m})|\!: \\[.2pc]
&\quad\ \ (u_{1},\dots,u_{m}),(x_{1},\dots,x_{m})\in I^{m}\\[.2pc]
&\quad\ \ \hbox{and }\left| u_{i}-x_{i}\right| \leq \delta
_{i},\text{ } i=1,2,\dots,m\}.
\end{align*}
Then let $H_{w_{m}}$ be the space of all real-valued functions $f$
satisfying
\begin{align*}
&\left| f(u_{1},\dots,u_{m})-f(x_{1},\dots,x_{m})\right|\\[.2pc]
&\quad\, \leq w_{2}\left(f;\left|
\frac{u_{1}}{1+u_{1}}-\frac{x_{1}}{1+x_{1}}\right|,\dots,\left|
\frac{u_{m}}{1+u_{m}}-\frac{x_{m}}{1+x_{m}}\right| \right).
\end{align*}
Therefore, using a similar technique in the proof of Theorem~2.1
one can obtain the following result immediately.

\begin{theor}[\!]
Let $A=(a_{jn})$ be a non-negative regular summability matrix{\rm
,} and let $ \{L_{n}\}$ be a sequence of positive linear operators
from $H_{w_{m}}$ into $ C_{B}(I^{m}).$ Then{\rm,} for any function
$f\in H_{w_{m}},$
\begin{equation*}
st_{A}-\lim_{n}\left\| L_{n}f-f\right\| =0
\end{equation*}
is satisfied if the following holds{\rm :}
\begin{equation*}
st_{A}-\lim_{n}\left\| L_{n}f_{i}-f_{i}\right\| =0,\quad
i=0,1,\dots,m+1,
\end{equation*}
where
\begin{align*}
&f_{0}(u_{1},\dots,u_{m}) = 1, f_{i}(u_{1},\dots,u_{m})=
\frac{u_{i}}{1+u_{i}},\quad i=1,2,\dots,m,\\[.2pc]
&f_{m+1}(u_{1},\dots,u_{m}) =
\sum_{k=1}^{m}\left(\frac{u_{k}}{1+u_{k}} \right) ^{2}.
\end{align*}
\end{theor}

\begin{remar}
If we choose $m=1,$ $I=[0,\infty)$ and replace $ A=(a_{jn})$ by
the identity matrix, then Theorem~2.2 reduces to Theorem~A.
\end{remar}

\section{Concluding remarks}

\setcounter{equation}{0}

In this section we display an example such that Theorem~A does not
work but Theorem~2.1 does.

Let $A=(a_{jn})$ be a non-negative regular summability matrix such
that $ \lim_{j}\max_{k} \{a_{jn}\}=0.$ In this case we know from
\cite{Kol} that $A$-statistical convergence is stronger than
ordinary convergence. So, we can choose a non-negative sequence
$(u_{n})$ which converges $A$-statistically to $1$ but
non-convergent. Assume that $I=[0,\infty)$ and $K:=I^{2}=I\times
I$. We now consider the following positive linear operators
defined on $H_{w_{2}}(K)$:
\begin{align*}
\hskip -4pc
T_{n}(f;x,y)=\frac{u_{n}}{(1+x)^{n}(1+y)^{n}}\sum_{k=0}^{n}
\sum_{l=0}^{n}f \left(\frac{k}{n-k+1}, \frac{l}{n-l+1}\right)
\binom{n}{k}\binom{n}{l} x^{k}y^{l},
\end{align*}
where $f\in H_{w_{2}},$ $(x,y)\in K$ and $n\in \Bbb{N}$.

We should remark that in the case of $u_{n}=1$, the operators
$T_{n}$ turn out to be the operators of Bleimann, Butzer and Hahn
\cite{Ble} (of two variables).

Since
\begin{equation*}
(1+x)^{n}=\sum_{k=0}^{n}\binom{n}{k}x^{k},
\end{equation*}
it is clear that, for all $n\in \Bbb{N}$,
\begin{equation*}
T_{n}(f_{0};x,y)=u_{n}.
\end{equation*}
Now, by assumption we have
\begin{equation}
st_{A}-\lim_{n}\left\| T_{n}f_{0}-f_{0}\right\| =st_{A}-\lim_{n}
\left| u_{n}-1\right| =0.
\end{equation}
Using the definition of $T_{n},$ we get
\begin{align*}
T_{n}(f_{1};x,y) &=
\frac{u_{n}}{(1+x)^{n}(1+y)^{n}}\sum_{k=1}^{n}\frac{k}{
n+1}\binom{n}{k}x^{k}\sum_{l=0}^{n}\binom{n}{l}y^{l} \\[.2pc]
&= \frac{u_{n}x}{(1+x)^{n}}\sum_{k=0}^{n-1}
\frac{k+1}{n+1}\binom{n}{k+1} x^{k}.
\end{align*}
Since
\begin{equation*}
\binom{n}{k+1}=\frac{n}{k+1}\binom{n-1}{k},
\end{equation*}
we obtain
\begin{align*}
T_{n}(f_{1};x,y) &=
\frac{u_{n}x}{1+x}\frac{1}{(1+x)^{n-1}}\sum_{k=0}^{n-1}
\frac{n}{n+1}\binom{n-1}{k}x^{k} \\[.2pc]
&= \frac{nu_{n}}{n+1}\left(\frac{x}{1+x}\right),
\end{align*}
which yields
\begin{equation*}
\left| T_{n}(f_{1};x,y)-f_{1}(x,y)\right| =\frac{x}{1+x}\left|
\frac{n}{n+1} u_{n}-1\right|
\end{equation*}
and hence
\begin{equation}
\left\| T_{n}f_{1}-f_{1}\right\| \leq \left|
\frac{n}{n+1}u_{n}-1\right|.
\end{equation}
Since $\lim_{n}\frac{n}{n+1}=1$ and $st_{A}-\lim_{n}u_{n}=1,$
observe that $ st_{A}-\lim_{n}$ $\frac{n}{n+1}u_{n}=1,$ so it
follows from (3.2) that
\begin{equation}
st_{A}-\lim_{n}\left\| T_{n}f_{1}-f_{1}\right\| =0.
\end{equation}
Similarly, we get
\begin{equation}
st_{A}-\lim_{n}\left\| T_{n}f_{2}-f_{2}\right\| =0.
\end{equation}
Finally, we claim that
\begin{equation}
st_{A}-\lim_{n}\left\| T_{n}f_{3}-f_{3}\right\| =0.
\end{equation}
Indeed, by the definition of $T_{n}$ we have
\begin{align*}
\hskip -4pc T_{n}(f_{3}; x,y) &=
\frac{u_{n}}{(1+x)^{n}(1+y)^{n}}\sum_{k=0}^{n}
\sum_{l=0}^{n}\left[ \frac{k^{2}}{\left(n+1\right)
^{2}}+\frac{l^{2}}{\left(n+1\right) ^{2}}\right]
\binom{n}{k}\binom{n}{l}x^{k}y^{l}\\[.2pc]
&= \frac{u_{n}}{(1+x)^{n}(1+y)^{n}}\sum_{k=1}^{n}
\frac{k^{2}}{\left(n+1\right)^{2}}\binom{n}{k}x^{k}
\sum_{l=0}^{n}\binom{n}{l}y^{l}\\[.2pc]
&\quad\, +\frac{u_{n}}{(1+x)^{n}(1+y)^{n}}
\sum_{k=0}^{n}\binom{n}{k} x^{k}\sum_{l=1}^{\infty
}\frac{l^{2}}{\left(n+1\right)^{2}}\binom{n}{l}y^{l}\\[.2pc]
&=
\frac{u_{n}}{(1+x)^{n}}\sum_{k=2}^{n}\frac{k(k-1)}{\left(n+1\right)
^{2}} \binom{n}{k}x^{k}+
\frac{u_{n}}{(1+x)^{n}}\sum_{k=1}^{n}\frac{k}{\left(n+1\right) ^{2}}x^{k} \\[.2pc]
&\quad\, +\frac{u_{n}}{(1+y)^{n}}
\sum_{l=2}^{n}\frac{l(l-1)}{\left(n+1\right) ^{2}}
\binom{n}{l}y^{l}+\frac{u_{n}}{(1+y)^{n}}\sum_{l=1}^{n}\frac{l}{\left(n+1\right) ^{2}}y^{l}%\\[.2pc]
\end{align*}
\begin{align*}
&= \frac{u_{n}x^{2}}{(1+x)^{n}}
\sum_{k=2}^{n-2}\frac{(k+2)(k+1)}{\left(n+1\right)
^{2}}\binom{n}{k+2}x^{k}+\frac{u_{n}x}{(1+x)^{n}}\sum_{k=0}^{n-1}
\frac{k+1}{\left(n+1\right) ^{2}}x^{k}\\[.2pc]
&\quad\, +\frac{u_{n}y^{2}}{(1+y)^{n}}
\sum_{l=0}^{n-2}\frac{(l+2)(l+1)}{\left(n+1\right)
^{2}}\binom{n}{l+2}y^{l}+ \frac{u_{n}y}{(1+y)^{n}}\sum_{l=1}^{n-1}
\frac{l+1}{\left(n+1\right) ^{2}}y^{l}.
\end{align*}
Now using the facts that
\begin{align*}
\binom{n}{k+2} &= \frac{n(n-1)}{(k+1)(k+2)}\binom{n-2}{k},
\binom{n}{k+1}=\frac{n}{k+1}\binom{n-1}{k}, \\[.2pc]
\binom{n}{l+2} &= \frac{n(n-1)}{(l+1)(l+2)}\binom{n-2}{l},
\binom{n}{l+1}=\frac{n}{l+1}\binom{n-1}{l},
\end{align*}
we get
\begin{align*}
T_{n}(f_{3};x,y) &=
\frac{n(n-1)u_{n}}{(n+1)^{2}}\frac{x^{2}}{(1+x)^{2}}
\frac{1}{(1+x)^{n-2}}\sum_{k=0}^{n-2}\binom{n-2}{k}x^{k} \\[.2pc]
&\quad\,
+\frac{nu_{n}}{(n+1)^{2}}\frac{x}{1+x}\frac{1}{(1+x)^{n-1}}\sum_{k=0}^{n-1}
\binom{n-1}{k}x^{k}\\[.2pc]
&\quad\, \times
\frac{n(n-1)u_{n}}{(n+1)^{2}}\frac{y^{2}}{(1+y)^{2}}\frac{1}{(1+y)^{n-2}}
\sum_{k=0}^{n-2}\binom{n-2}{l}y^{l} \\[.2pc]
&\quad\,
+\frac{nu_{n}}{(n+1)^{2}}\frac{y}{1+y}\frac{1}{(1+y)^{n-1}}\sum_{k=0}^{n-1}
\binom{n-1}{l}y^{l}\\[.2pc]
&=
\frac{n(n-1)u_{n}}{(n+1)^{2}}\frac{x^{2}}{(1+x)^{2}}+\frac{nu_{n}}{
(n+1)^{2}}\frac{x}{1+x} \\[.2pc]
&\quad\,
+\frac{n(n-1)u_{n}}{(n+1)^{2}}\frac{y^{2}}{(1+y)^{2}}+\frac{nu_{n}}{
(n+1)^{2}}\frac{y}{1+y},
\end{align*}
which implies that
\begin{align*}
\left| T_{n}(f_{3};x,y)-f_{3}(x,y)\right| &=
\left(\frac{x^{2}}{(1+x)^{2}}+ \frac{y^{2}}{(1+y)^{2}}\right)
\left|
\frac{n(n-1)u_{n}}{(n+1)^{2}}-1\right|\\[.2pc]
&\quad\, +\left(\frac{x}{1+x}+\frac{y}{1+y}\right)
\frac{nu_{n}}{(n+1)^{2}}.
\end{align*}
Taking supremum over $(x,y)\in K$ we have
\begin{equation}
\left\| T_{n}f_{3}-f_{3}\right\| \leq 2(\alpha _{n}+\beta _{n}),
\end{equation}
where $\alpha _{n}:=\big| \frac{n(n-1)u_{n}}{(n+1)^{2}}-1 \big|$
and $ \beta _{n}:=\frac{nu_{n}}{(n+1)^{2}}.$ Then, by assumption,
observe that
\begin{equation}
st_{A}-\lim_{n}\alpha _{n}=st_{A}-\lim_{n}\beta _{n}=0.
\end{equation}
Now given $\varepsilon >0$, define the following sets
\begin{align*}
U &:= \{n\hbox{:}\ \| T_{n}f_{3}-f_{3}\| \geq \varepsilon \},\\[.2pc]
U_{1} &:=\left\{ n\hbox{:}\ \alpha _{n}\geq \frac{\varepsilon
}{4}\right\},\quad U_{2}:= \left\{n\hbox{:}\ \beta _{n}\geq
\frac{\varepsilon }{4}\right\}.
\end{align*}
By (3.6), it is obvious that $U\subseteq U_{1}\cup U_{2}.$ Then,
for each $ j\in \Bbb{N}$, we may write that
\begin{equation}
\sum_{n\in U}a_{jn}\leq \sum_{n\in U_{1}}a_{jn}+\sum_{n\in
U_{2}}a_{jn}.
\end{equation}
Taking limit as $j\rightarrow \infty $ in (3.8) and using (3.7) we
have
\begin{equation*}
\lim_{j}\sum_{n\in U}a_{jn}=0,
\end{equation*}
which proves (3.5).

Therefore, using (3.1), (3.3), (3.4) and (3.5) in Theorem~2.1, we
obtain that, for all $f\in H_{w_{2}},$
\begin{equation*}
st_{A}-\lim_{n}\| T_{n}f-f\| =0.
\end{equation*}
However, since the sequence $(u_{n})$ is non-convergent, $\{T_{n}f\}$ is not
uniformly convergent to $f.$


\begin{thebibliography}{99}
\bibitem{Alt} Altomare F and Campiti M, Korovkin type approximation theory
and its applications, {\it de Gruyter Stud. Math.} (Berlin: de
Gruyter) (1994) vol.~17

\bibitem{Ble} Bleimann G, Butzer P L and Hahn L, A Bernstein type operator
approximating continuous functions on semiaxis, {\it Indag. Math.}
{\bf 42} (1980) 255--262

\bibitem{Boj} Bojanic R and Khan M K, Summability of Hermite--Fej\'{e}r
interpolation for functions of bounded variation, {\it J. Nat.
Sci. Math.} {\bf 32} (1992) 5--10

\bibitem{Boos} Boos J, Classical and modern methods in summability (UK: Oxford
University Press) (2000)

\bibitem{Che} Cheney E W and Sharma A, Bernstein power series, {\it
Canad. J. Math.} {\bf 16} (1964) 241--253

\bibitem{C-G} \c{C}akar \"{O} and Gadjiev A D, On uniform approximation
by Bleimann, Butzer and Hahn on all positive semiaxis, {\it Trans.
Acad. Sci. Azerb. Ser. Phys. Tech. Math. Sci.} {\bf 19} (1999)
21--26

\bibitem{Dum} Duman O, Khan M K and Orhan C, $A$-statistical convergence of
approximating operators, {\it Math. Inequal. Appl.} {\bf 6} (2003)
689--699

\bibitem{Fas} Fast H, Sur la convergence statistique, {\it Colloq. Math.}
{\bf 2} (1951) 241--244

\bibitem{Fre} Freedman A R and Sember J J, Densities and summability,
{\it Pacific J. Math.} {\bf 95} (1981) 293--305

\bibitem{Fri} Fridy J A, On statistical convergence, {\it Analysis}
{\bf 5} (1985) 301--313

\bibitem{Gad} Gadjiev A D and Orhan C, Some approximation theorems via
statistical convergence, {\it Rocky Mountain J. Math.} {\bf 32}
(2002) 129--138

\bibitem{Khan} Khan M K, Vecchia B D and Fassih A, On the monotonicity of
positive linear operators, {\it J. Approx. Theory} {\bf 92} (1998)
22--37

\bibitem{Kol} Kolk E, Matrix summability of statistically convergent
sequences, {\it Analysis} {\bf 13} (1993) 77--83

\bibitem{Kol-2} Kolk E, The statistical convergence in Banach spaces,
{\it Tartu \"{U}l. Toimetised} {\bf 928} (1991) 41--52

\bibitem{Kor} Korovkin P P, Linear operators and approximation theory
(Delhi: Hindustan Publ. Co.) (1960)

\bibitem{Mil} Miller H I, A measure theoretical subsequence
characterization of statistical convergence, {\it Trans. Am. Math.
Soc.} {\bf 347} (1995) 1811--1819

\bibitem{Ste} Steinhaus H, Sur la convergence ordinaire et la
convergence asymptotique, {\it Colloq. Math.} {\bf 2} (1951)
73--74
\end{thebibliography}
\end{document}